\documentclass[10pt,reqno]{amsart}

\usepackage{fixltx2e}                       
\usepackage[usenames,dvipsnames]{xcolor}
\usepackage{fancyhdr}
\usepackage{amsmath,amsfonts,amsbsy,amsgen,amscd,amssymb,amsthm}
\usepackage{mathrsfs,eucal}
\usepackage{mathtools}
\usepackage[font=small,margin=0.25in,labelfont={sc},labelsep={colon}]{caption}
\usepackage{tikz}
\usepackage{enumitem}
\usepackage[colorlinks=true,linkcolor=blue,citecolor=green]{hyperref}
\usepackage{graphicx}

\usepackage[T1]{fontenc}
\usepackage{bm}

\newtheorem{theorem}{Theorem}[section]

\newtheorem{definition}[theorem]{Definition}
\newtheorem{remark}[theorem]{Remark}
\newtheorem{algorithm}[theorem]{Algorithm}
\newtheorem{proposition}[theorem]{Proposition}

\DeclareMathOperator{\Tr}{\mathrm{Tr}}

\DeclareMathOperator{\GGG}{\mathcal{G}}
\DeclareMathOperator{\rank}{\operatorname{rank}}
\DeclareMathOperator{\1}{\mathbf{1}}

\DeclareMathOperator{\diag}{\mathrm{diag}}

\newcommand{\defeq}{\mathrel{\mathop:}=}

\begin{document}

%\title[Probably Certifiably Correct Algorithms]{Probably Certifiably Correct Algorithms}
\title{A note on Probably Certifiably Correct algorithms}
%\title{Probably Certifiably Correct Algorithms}
\author{Afonso S. Bandeira}
\address[Bandeira]{Department of Mathematics, Massachusetts Institute of Technology, Cambridge, Masachusetts 02142, USA ({\tt bandeira@mit.edu}).}
\thanks{A.~S.~Bandeira was partially supported by AFOSR Grant No.\ FA9550-12-1-0317 and NSF grant DMS-1317308. Part of this work was done while the author was at Princeton University.}
%\thanks{ASB was supported by AFOSR Grant No.\ FA9550-12-1-0317.}
%\email{ajsb@math.princeton.edu}

\begin{abstract}
Many optimization problems of interest are known to be intractable, and while there are often heuristics that are known to work on typical instances, it is usually not easy to determine a posteriori whether the optimal solution was found. In this short note, we discuss algorithms that not only solve the problem on typical instances, but also provide a posteriori certificates of optimality, probably certifiably correct (PCC) algorithms. As an illustrative example, we present a fast PCC algorithm for minimum bisection under the stochastic block model and briefly discuss other examples.
\end{abstract}

%\subjclass[2000]{Primary 60B20; Secondary 46B09, 60F10}


\maketitle

\thispagestyle{empty}

%\tableofcontents

\section{Introduction}

Estimation problems in many areas are often formulated as optimization problems, maximum likelihood estimation being a prime example. Unfortunately, these optimization problems are, in many relevant instances, believed to be computationally intractable.

To circumvent the computational intractability of these problems, an impressive body of work is dedicated to propose and understand algorithms that are guaranteed to produce a solution with certain approximation guarantees~\cite{DPWilliamsonDBShmoys_2011}. %While these guarantees hold for any input, they are usually formulated in terms of the objective value of the solution and so not always directly translate to meaningful statement about the statistical properties of the solution. Also, 
Since these guarantees have to hold for worst-case inputs, they tend to be pessimistic and, fortunately, many real-world scenarios do not resemble these worst-case instances.

In line with the ``real-world data is not your enemy'' paradigm, there is a line of work that attempts to propose and understand algorithms that work on only some sets of, hopefully typical, inputs. A prime example is sparse recovery, where the popular Compressed Sensing papers of Cand\`{e}s, Romberg, Tao, and Donoho~\cite{Candes_CS1,Donoho_CS} established that, while finding a sparse solution to an underdetermined system is computationally intractable in the worst-case, a simple efficient algorithm succeeds with high probability in many natural probabilistic models of instances. Important examples also include planted bisection~\cite{McSherry,Feige_Kilian_bisection_01} and matrix completion~\cite{candes2009exact}, among many others.

\begin{definition}[Probabilistic Algorithm]
Given an optimization problem that depends on an input and a probability distribution $D$ over the inputs, we say an algorithm is a probabilistic algorithm for this problem and distribution of instances if it finds an optimal solution with high probability\footnote{We say that an event happens with high probability if its probability tends to $1$ as the underlying parameters tend to infinity.}, with respect to $D$.
\end{definition}

Although these guarantees make excellent arguments towards the efficacy of certain probabilistic algorithms, they have the drawback that, oftentimes, it is not easy to check a posteriori whether the solution computed is the optimal one. Even in the idealized scenario were the input really follows the distribution in the guarantee there is a nonzero probability of the algorithm not producing the optimal solution, and it is often not possible to tell whether it did or not. To make matters worse, in practice, the exact distribution of instances rarely exactly matches the idealized one in the guarantees.

The situation is different for a certain class of algorithms, convex relaxation based ones~\cite{CandesTao_ConvexMatCompletion10}. Some of these methods work by enlarging the feasibility set of the optimization problem to a convex set where optimizing the objective function becomes tractable. While the optimal solution is not guaranteed to be in the original feasibility set, there are many examples for which rounding procedures are known to produce solutions with approximation guarantees (one such example being the Goemans-Williamson~\cite{MXGoemans_DPWilliamson_1995} approximation algorithm for \texttt{Max-Cut}). On the other hand, if the solution happens to lie on the original feasibility set, then one is sure that it must be the optimal solution of the original problem (providing, also, an a posteriori certificate). Fortunately, this tends to be the case for many examples of problems and relaxations~\cite{Amini_Wainwright_SparsePCA,Awasthi_RelaxnoRound,Bandeira_rankrecoveryangsynch,Bandeira2014open,Candes_Strohmer_Voroninski_phaselift}. This motivates the following definition.

\begin{definition}[Probably Certifiably Correct (PCC) Algorithm]
Given an optimization problem that depends on an input and a probability distribution $D$ over the inputs, we say an algorithm is a Probably Certifiably Correct (PCC) algorithm for this problem and distribution of instances if, with high probability (w.r.t. $D$), it finds an optimal solution and certifies to have done so. Moreover, it never incorrectly certifies a non-optimal solution.
\end{definition}

A PCC algorithm has the advantage of being able to produce a posteriori certificates. In particular, this renders them more appealing to be used in examples where the distribution of problem instances may not coincide with the idealized ones in proved guarantees.
Indeed, duality is very often used in practice to provide a posteriori certificates of quality of a solution to an optimization problem (see~\cite{Rosen_SLAMduality} for a particularly recent example in the problem of Simultaneous Localization and Mapping (SLAM)).
While this is a great argument towards the use of convex relaxation-based approaches,\footnote{Note that there is another type of convex relaxation-based algorithms, such as in sparse recovery, where instead of enlarging the feasibility set, one replaces the objective function by a complex surrogate. Unfortunately, in that case, it appears to be more difficult to certify, a posteriori, the optimality of the solutions. For the particular case of sparse recovery, we refer the reader to~\cite{Gribonval_CertificateSparseRecovery} for a discussion on certificates of optimality.} the convex problems that many such algorithms are relaxed to are semidefinite programs (SDP) which, while solvable (to arbitrary precision) in polynomial time~\cite{LVanderberghe_SBoyd_1996}, tend to be computationally costly.\footnote{It is also fairly common for problems to be relaxed to linear programs, which tend to be computationally cheaper.}

Many probabilistic algorithms not based on convex relaxations (such as, for example, spectral methods~\cite{McSherry}, message-passing type algorithms~\cite{Montanari_positivePCA_iter}, alternating-minimization/expectation-maximization techniques~\cite{Jain:matrixcompletion_altmin_13}) often do not enjoy a posteriori guarantees, but tend to be considerably more efficient than the convex relaxation-based competitors. This motivates the natural question of whether it is possible to devise a posteriori certifiers for candidate solutions produced by these or other methods.

\begin{definition}[Probabilistic Certifier]
Given an optimization problem that depends on an input, a probability distribution $D$ over the inputs, and a candidate solution for it, we call a Probabilistic Certifier, a method that:
\begin{itemize}
 \item With high probability (w.r.t. $D$), if the candidate solution is an optimal one it outputs: \texttt{The solution is optimal}.\footnote{In some cases, certifiers may also certify that a solution is not only optimal, but the unique optimal solution, as it will be the case with Algorithm~\ref{alg:PCCSBMcheckPSD}.} It may, with vanishing probability,\footnote{By vanishing probability we mean probability tending to $0$ as the underlying parameters tend to infinity.} output: \texttt{Not sure whether the solution is optimal}. 
 \item If the candidate solution is not an optimal solution, it always outputs: \texttt{Not sure whether the solution is optimal}.
\end{itemize}
\end{definition}

A particularly natural way of constructing such certifiers is to rely on convex relaxation-based PCC algorithms; given a candidate solution computed by a probabilistic algorithm, one can check whether it is an optimal solution to a certain convex relaxation. Remarkably, it is sometimes considerably faster to check wether a candidate solution is an optimal solution of a convex program than to solve the program; in many such cases, one can devise faster PCC algorithms by combining fast probabilistic algorithms with fast methods to certify that a candidate solution is an optimal solution to a convex relaxation, or even that it is the unique optimal solution (as it will be the case with Algorithm~\ref{alg:PCCSBMcheckPSD}).
In the next section, we use the the problem of minimum bisection under the stochastic block model to illustrate these ideas.

\section{A fast PCC algorithm for recovery in the Stochastic Block Model}

The problem of minimum bisection on a graph is a particularly simple instance of community detection that is known to be NP-hard. Recently, there has been interest in understanding the performance of several heuristics in typical realizations of a certain random graph model that exhibits community structure, the stochastic block model: given $n$ even and $0\leq q < p \leq 1$, we say that a random graph $G$ is drawn from $\GGG(n,p,q)$, the Stochastic Block Model with two communities, if $G$ has $n$ nodes, divided in two clusters of $\frac{n}2$ nodes each, and for each pair of vertices $i,j$, $(i,j)$ is an edge of $G$ with probability $p$ if $i$ and $j$ are in the same cluster and with probability $q$ otherwise, independently from any other edge. 

Let $A$ denote the adjacency matrix of $G$. We define the signed adjacency matrix $B$ as $B = 2A - (\1\1^T - I)$. To each partitioning of the nodes we associate a vector $x$ with $\pm1$ entries corresponding to cluster memberships. The minimum bisection of $G$ can be written as
\begin{equation}\label{eq:min_bissection}
 \begin{array}{rl}
  \displaystyle{\max_x} & \displaystyle{x^T B x} \\ % \displaystyle{\sum_{i,j=1}^{n}B_{ij}x_ix_j} \\ 
  \text{ s.t. } & %x \in \RR^n\\
  %&
   \displaystyle{x_i^2 = 1,\ \forall_i}\\
  & \displaystyle{\sum_{i=1}^nx_i = 0.}
 \end{array}
\end{equation}

Setting $p = \alpha\frac{\log n}n$ and $q = \beta \frac{\log n}n$, it is known~\cite{Abbe_SBMExact,Mossel_SBM3_exact} that the hidden partition, with high probability, coincides with the minimum bisection, and can be computed efficiently provided that
\begin{equation}\label{eq:regime_SBM_possible}
 \sqrt{\alpha} - \sqrt{\beta} > \sqrt{2}.
\end{equation}
On the other hand, if
\(
 \sqrt{\alpha} - \sqrt{\beta} < \sqrt{2},
\)
then, with high probability, the maximum likelihood estimator (which corresponds to the minimum bisection) does not coincide with the hidden partition.%\footnote{There is a fascinating body of literature on the problem of detection of the community structure, while not in the scope of this note, we refer the reader to~\cite{Massoulie_SBM,Mossel_SBM2} and references within for more on this problem.}

Remarkably, for the stochastic block model on two communities with parameters in the regime~\eqref{eq:regime_SBM_possible}, there are quasi-linear time algorithms known to be probabilistic algorithms for minimum bisection (see, for example,~\cite{AbbeSandon15}). A convex relaxation, proposed in~\cite{Abbe_SBMExact}, was also recently shown to exactly compute the minimum bisection in the same regime~\cite{Bandeira_Laplacian,Hajek_et_al_SBM_SDP}.

The convex relaxation in~\cite{Abbe_SBMExact,Bandeira_Laplacian,Hajek_et_al_SBM_SDP} is obtained by writing~\eqref{eq:min_bissection} in terms of a new variable $X=xx^T$. More precisely,~\eqref{eq:min_bissection} is equivalent to
\begin{equation}\label{SBM1_lifted}
 \begin{array}{rl}
  \displaystyle{\max_X} & \displaystyle{\Tr\left(BX\right)} \\ 
  \text{ s.t. } & \displaystyle{X_{ii} = 1,\ \forall_i} \\
    & \displaystyle{X \succeq 0} \\
    & \displaystyle{\rank(X) = 1} \\
    & \displaystyle{\Tr\left(X\1\1^T\right) = 0}.
 \end{array}
\end{equation}
The semidefinite programming relaxation considered is obtained by removing the last two constraints.~\footnote{\eqref{SBM1_SDP} would still be a semidefinite program if the constraint $\Tr\left(X\1\1^T\right) = 0$ was kept, but it turns out that it is not needed and removing it renders the analysis slightly simpler.}
\begin{equation}\label{SBM1_SDP}
 \begin{array}{rl}
  \displaystyle{\max_X} & \displaystyle{\Tr\left(BX\right)} \\ 
  \text{ s.t. } & \displaystyle{X_{ii} = 1,\ \forall_i} \\
    & \displaystyle{X \succeq 0}.
 \end{array}
\end{equation}
The argument in~\cite{Bandeira_Laplacian,Hajek_et_al_SBM_SDP} makes use of duality~\cite{LVanderberghe_SBoyd_book}. Since~\eqref{SBM1_SDP} satisfies Slater's condition,\footnote{Slater's condition is a technical condition that ensures strong duality and tends to be satisfied in many relevant problems; in this case it asks that there is a feasible point $X$ that has no zero eigenvalues and the identity matrix serves as an example, see~\cite{LVanderberghe_SBoyd_book} for more details.} the optimal value of dual program given by
\begin{equation}\label{SBM1_SDP_dual}
 \begin{array}{rl}
  \displaystyle{\min_D} & \displaystyle{\Tr(D)} \\ 
  \text{ s.t. } & \displaystyle{D - B \succeq 0} \\
   & \displaystyle{D \text{ is diagonal}}
 \end{array}
\end{equation}
is known to match the optimal value of~\eqref{SBM1_SDP}. More precisely, given the hidden partition $x_\natural\in\{\pm1\}^{n}$, Abbe et al.~\cite{Abbe_SBMExact} propose $D_\natural \defeq D_{\diag(x_\natural)B\diag(x_\natural)}$ as a candidate solution for the dual,\footnote{Interestingly, in this case, the equality constraints and complementary slackness conditions are enough to pin-point a single possible candidate for a dual certificate, and only the semidefinite constraint needs to be checked; this is the case for semidefinite programs satisfying certain properties, see~\cite{alizadeh1997complementarity} for more details.} where $\diag(x_\natural)$ is a diagonal matrix whose diagonal is given by $x_\natural$ and $D_\natural = D_{\diag(x_\natural)B\diag(x_\natural)}$ is a diagonal matrix whose diagonal is given by
\[
 \left[D_\natural\right]_{ii} = \left[D_{\diag(x_\natural)B\diag(x_\natural)}\right]_{ii} = \sum_{j=1}^n \left[\diag(x_\natural)B\diag(x_\natural)\right]_{ij} = \left(x_\natural\right)_i\sum_{j=1}^n B_{ij}\left(x_\natural\right)_j.
\]
More recently,~\cite{Bandeira_Laplacian,Hajek_et_al_SBM_SDP} showed that, in the parameter regime given by~\eqref{eq:regime_SBM_possible} and with high probability, this dual candidate solution is indeed a feasible solution to~\eqref{SBM1_SDP_dual} whose value matches the value of the $x_\natural^TBx_\natural$. This implies that $x_\natural x_\natural^T$ is an optimal solution of \eqref{SBM1_SDP}. Moreover, since~\cite{Bandeira_Laplacian,Hajek_et_al_SBM_SDP} show that
\begin{equation}\label{conditionsPSD_SBM}
D_{\natural} - B\succeq 0 \quad \text{ and } \quad \lambda_2\left(D_{\natural}- B\right) > 0,
\end{equation}
where $\lambda_2$ denotes the second smallest eigenvalue, the argument can be easily strengthened (using complementary slackness) to show that $x_\natural x_\natural^T$ is the unique optimal solution (see~\cite{Bandeira_Laplacian,Hajek_et_al_SBM_SDP} for details).

A particularly enlightening way of showing that~\eqref{conditionsPSD_SBM} indeed certifies that the partitioning given by $x_\natural$ is the unique solution to the minimum bisection problem is to note that, for any other candidate bisection $x\in\{\pm 1\}^n$,
\begin{eqnarray}
 x_\natural^T B x_\natural - x^T B x & = & x^T \left[ D_{\natural} - B  \right]^Tx + \sum_{i=1}^n\left[ D_{\natural}  \right]_{ii}\left(1-x_i^2\right) \nonumber\\
  & = & x^T \left[ D_{\natural} - B  \right]x. \label{dual_certificate_SOS}
\end{eqnarray}

Since $D_{\natural} - B$ is known to satisfy~\eqref{conditionsPSD_SBM}, then $x^T \left[ D_{\natural} - B  \right]x \geq 0$. Moreover, since $\left[ D_{\natural} - B  \right]x_\natural = 0$, if $x$ corresponds to another bisection (meaning that $x\neq x_\natural$ and $x\neq -x_\natural$) then $x^T \left[ D_{\natural} - B  \right]x > 0$, implying that $x_\natural^T B x_\natural - x^T B x > 0$.

\begin{remark}
 A particularly fruitful interpretation of~\eqref{dual_certificate_SOS} is to think of it as a sum-of-squares certificate. More precisely, since $D_{\natural} - B$ is positive semidefinite it has a Cholesky decomposition $D_{\natural} - B = VV^T$ which means that, for $x\in\{\pm 1\}^n$,
 \[
% x^T \left[ D_{\natural} - B  \right]x
 x_\natural^T B x_\natural - x^T B x = x^TVV^Tx = \left\|V^Tx\right\|^2 = \sum_{j=1}^n \left(\sum_{i=1}^n V_{ij}x_i\right)^2.
 \]
 By writing $x_\natural^T B x_\natural - x^T B x$ has a sum of squares, we certify that $x_\natural^T B x_\natural$ is an optimal solution. It turns out that certificates of this type always exist, potentially having to include polynomials of larger degree, and that, with a fixed bound on the degree of the polynomials involved, these certificates can be found with semidefinite programming whose complexity depends on the degree bound. This is a simple instance of the sum-of-squares technique (based on Stengle's Positivstellensatz~\cite{Stengle_Positivstellensatz}) % in real algebraic geometry,
  proposed independently in a few different areas~\cite{Shor_87_SOS,Nesterov_00_SOS,Lassere_01_SOS,Parrilo_thesis_SOS} and now popular in theoretical computer science~\cite{Barak_Steurer_surveyICM}.
\end{remark}

This suggests the following PCC algorithm for minimum bisection in the Stochastic Block Model.

\begin{algorithm}\label{alg:PCCSBMcheckPSD}
 Given a graph $G$, use the quasi-linear time algorithm in~\cite{AbbeSandon15} to produce a candidate bisection $x_\ast$. If $x_\ast^T1=0$ and
 \begin{equation}\label{eq:conditionxast}
  \lambda_2\left( D_{\diag(x_\ast)B\diag(x_\ast)}  - B \right) > 0,
 \end{equation}
output:
\begin{itemize}
 \item \texttt{$x_\ast$ is the minimum bisection}.
\end{itemize}
If not, output:
\begin{itemize}
 \item \texttt{Not sure whether $x_\ast$ is the minimum bisection}.
\end{itemize}
\end{algorithm}

Note that, since $\left( D_{\diag(x_\ast)B\diag(x_\ast)}  - B \right)x_\ast = 0$, \eqref{eq:conditionxast} automatically implies condition~\eqref{conditionsPSD_SBM}.

The following follows immediately from the results in~\cite{AbbeSandon15} and~\cite{Bandeira_Laplacian,Hajek_et_al_SBM_SDP}.

\begin{proposition}
Algorithm~\ref{alg:PCCSBMcheckPSD} is a Probably Certifiably Correct Algorithm for minimum bisection under the stochastic block model in the regime of parameters given by~\eqref{eq:regime_SBM_possible}.
\end{proposition}

\subsection{Randomized certificates}

While Algorithm~\ref{alg:PCCSBMcheckPSD} is considerably faster than solving the semidefinite program~\eqref{SBM1_SDP} it still requires one to check that an $n\times n$ matrix has a positive second-smallest eigenvalue, which we do not know how to do in quasi-linear time. A potentially faster alternative would be to use a randomized power-method-like algorithm, such as randomized Lanczos method~\cite{Kuczynski_PowerMethodEigenvalue}, to estimate the second smallest eigenvalue of $D_{\diag(x_\ast)B\diag(x_\ast)}  - B$. Note that since $B = 2A - (\1\1^T - I)$, where $A$ is a sparse matrix, matrix-vector multiplies with $D_{\diag(x_\ast)B\diag(x_\ast)}  - B$ can be computed in quasi-linear time. While the use of such randomized methods would not provide a probabilistic certificate, it could potentially provide a randomized certificate that has a small probability (with respect to a source of randomness independent to $D$) of ``certifying'' an incorrect solution. However, since it would rely on independent randomness, the process would be able to be repeated to achieve an arbitrarily small probability of providing false certificates.

% \begin{definition}[Randomized Certificate]
%  Given an optimization problem that depends on an input and a probability distribution $D$ over the inputs and a candidate solution for it, we refer to as a randomized Probabilistic Certifier, a method that has access to a source of randomness independent from the randomness of the problem instances and:
% \begin{itemize}
%  \item With high probability (w.r.t. $D$), and positive probability with respect to the independent randomness, if the candidate solution is the optimal one it outputs: \texttt{The solution is optimal}.
%  \item Whether the solution is optimal of not, with positive probability on the independent randomness it outputs: \textt{Am not sure whether I can certify if the s
%  It may, with vanishing probability, output: \texttt{Not sure whether the solution is optimal}. 
%  \item If the candidate solution is not the optimal solution, and with positive probability on the independent randomness, it always outputs: \texttt{Cannot certify that the solution is optimal.}.
% \end{itemize}
% \end{definition}

In fact, we believe that the analysis of the typical behavior of the second eigenvalue of $D_{\diag(x_\ast)B\diag(x_\ast)}  - B$ in~\cite{Bandeira_Laplacian,Hajek_et_al_SBM_SDP} and the guarantees for the performance of randomized Lanczos method~\cite{Kuczynski_PowerMethodEigenvalue} can be used to devise a quasi-linear time randomized procedure that can serve as a randomized certificate for minimum bisection in the stochastic block model, as described above. However, such construction falls outside of the scope of this short note, and so it is left for future research.

%$\OOO\left(\left(\frac{\log n}{k}\right)^2\right)$

% Even if we are given the dual variables $y_i$ we still have to verify whether $\sum_i y_i A_i - C \succeq 0$ and although much faster than solving an SDP it is not clear how to do it extremely fast. This poses an interesting question of certifying, fast, that a matrix is PSD.
% It is easy to do it with a randomized algorithm. Power Method (actually Lanczos algorithm) works well for this~\cite{Kuczynski_PowerMethodEigenvalue}, after k matrix multiplies you get an error of .
% 
% Randomized algorithm is different from probabilistic, it can be tried over and over again and one gets more sure that they worked

\section{Other examples of PCC algorithms and future directions}

One of the most appealing characteristics of Algorithm~\ref{alg:PCCSBMcheckPSD} above is that it can be easily generalized to many other settings. In fact, given an optimization problem and a distribution $D$ over the instances, if there is a fast probabilistic algorithm and a convex relaxation that is known to be tight (meaning that its optimal solution is feasible in the original problem), one can make use of both algorithms, similarly to Algorithm~\ref{alg:PCCSBMcheckPSD}, and devise a fast PCC algorithm: by first running the fast probabilistic algorithm and then checking whether the candidate solution is the optimal solution to the convex relaxation. Unfortunately, it is not clear, in general, whether one can check optimality in the convex relaxation considerably faster than simply solving it. On the other hand, many proofs of tightness of convex relaxations also provide a candidate dual solution and, in many instances, checking whether this dual solution is indeed a dual certificate is significantly faster.

For some problems, such as multisection in the stochastic block model with multiple communities, both fast probabilistic algorithms~\cite{AbbeSandon15} and tightness guarantees for convex relaxations have already been established~\cite{Hajek_et_al_SBM_SDP_extensions,Agawarl_multisection_SBM,PerryWein_MultisectionSBM} suggesting that this framework could be easily applied there. Other problems for which convex relaxations are known to be tight include: Synchronization over $\mathbb{Z}_2$~\cite{Abbe_Z2Synch} and $SO(2)$~\cite{Bandeira_rankrecoveryangsynch}, sparse PCA~\cite{Amini_Wainwright_SparsePCA}, k-medians and k-means clustering~\cite{Awasthi_RelaxnoRound,DustinSole_kmeansSDP1}, multiple-input multiple-output (MIMO) channel detection~\cite{AMCSo_MIMOSDP_10}, sensor network localization~\cite{AMCSo_Ye_SNL_06}, shape matching~\cite{Chen_Huang_Guibas_Graphics}, and many others.
There are several others where convex relaxations are conjectured to be tight under appropriate conditions, such as the non-negative PCA problem~\cite{Montanari_positivePCA_iter}, the multireference alignment problem~\cite{Bandeira_Charikar_Singer_Zhu_Alignment,Bandeira_NonUniqueGames}, and the synchronization problem in SLAM~\cite{Rosen_SLAMduality}. For the case of non-negative PCA, there is a known probabilistic algorithm~\cite{Montanari_positivePCA_iter}.
We suspect that this framework may be useful in devising fast PCC algorithms for a large class of problems, perhaps including many of the described above.\footnote{
An alternative strategy to developing fast PCC algorithms is to exploit particular structural properties of the problem to devise fast methods to solve the corresponding convex relaxations (see~\cite{Hopkins_Shi_Steurer_fastTensorPCASOS} for a recent example in tensor PCA).} Moreover, even when probabilistic algorithms are not available, fast certifiers may be useful to test the performance of heuristics in real world problems.

%A natural question raised by the observations on this note is whether there are good ways of constructing certificates that are not based on convex duality. In this direction, we refer to a certificate present in~\cite{Boppana_SBM1} that appears to be of a different nature than the ones described in this note, although it also relies in that although does not rely on convex duality quite in the same way although, in some sense, it also relies on convex optimization.

%Another example is Nicolas approach to solving SDPs~\cite{Boumal_RiemannianStairCase} \cite{Journee_staircased1}.

\subsection*{Acknowledgements}
The author thanks Dustin G. Mixon, Soledad Villar, Nicolas Boumal, and Amit Singer for interesting discussions and valuable comments on an earlier version of this manuscript. The author also acknowledges Dustin G. Mixon for suggesting the term \emph{Probably Certifiably Correct}. %The author presented some of these ideas  of these results in various seminars throughout the end of 2014 and beginning of 2015. Many questions and comments raised by the audience greatly improved the quality of this manuscript, a warm thanks to all of them. Thank Berkeley for the talk.

\bibliographystyle{siam}
\bibliography{../../afonso}

\end{document}